\documentclass[fleqn]{mat01}
\usepackage{times,mathtimy,amssymb,latexsym}%,krepsf}
\begin{document}

\setcounter{page}{347}
\firstpage{347}

\def\d{\hbox{d}}
\def\e{\hbox{e}}

\newtheorem{theor}{\bf Theorem}
\newtheorem{lem}{Lemma}
\newtheorem{rem}{Remark}
\newtheorem{pot}{\it Proof of Theorem}

\renewcommand{\theequation}{\thesection\arabic{equation}}

\title{Wavelet characterization of H\"ormander symbol class
$\pmb{S^{m}_{\rho,\delta}}$ and applications}

\markboth{Q~X~Yang}{Wavelet characterization of H\"ormander symbol
class $S^{m}_{\rho,\delta}$}

\address{Department of Mathematics, Wuhan University, 430072 Hubei, China\\
\noindent E-mail: yangqi99@public.wh.hb.cn}

\author{Q~X~YANG}

\volume{115}

\mon{August}

\parts{3}

\pubyear{2005}

\Date{MS received 12 March 2005}

\begin{abstract}
In this paper, we characterize the symbol in H\"ormander symbol
class $S^{m}_{\rho,\delta}\, (m\in R, \rho,\delta\geq 0)$ by its
wavelet coefficients. Consequently, we analyse the
kernel-distribution property for the symbol in the symbol class
$S^{m}_{\rho,\delta}\, (m\in R, \rho>0, \delta\geq 0)$  which is
more general than known results; for non-regular symbol operators,
we establish sharp $L^{2}$-continuity which is better than
Calder\'on and Vaillancourt's result, and establish $L^{p}\,
(1\leq p\leq\infty)$ continuity which is new and sharp. Our new
idea is to analyse the symbol operators in phase space with
relative wavelets, and to establish the kernel distribution
property and the operator's continuity on the basis of the wavelets
coefficients in phase space.
\end{abstract}


\maketitle

\setcounter{equation}{0}

\section{Introduction}

A symbol $\sigma(x,\xi)\in S'(R^{n}\times R^{n})$ can define a
symbol operator $\sigma(x,D)\hbox{:}\ S(R^{n})\rightarrow S'(R^{n})$ by
the following formula:
\begin{equation}
\sigma(x,D)f(x) = \int \e^{ix\xi}\sigma(x,\xi)\hat{f}(\xi) \ \d
\xi,
\end{equation}
where $\hat{f}(\xi)$ is the Fourier transformation of function
$f(x)$. When H\"ormander studied pseudodifferential operators, he
introduced H\"ormander's symbol class $S^{m}_{\rho,\delta}\, (m\in
R, \rho,\delta\geq 0)$. One writes $\sigma(x,\xi)\in
S^{m}_{\rho,\delta}$, if
\begin{equation}
|\partial^{\alpha}_{x}\partial^{\beta}_{\xi} \sigma(x,\xi)| \leq
C_{\alpha,\beta} (1+|\xi|)^{m-\rho |\beta| + \delta |\alpha|}, \ \
\forall \alpha,\beta\in N^{n}.
\end{equation}
But we did not know what are the elements in H\"ormander class
$S^{m}_{\rho,\delta}\, (m\in R, \rho,\delta\geq 0)$ before.
Professor Meyer \cite{12} proposed me to study such a kind of
pseudodifferential operators with wavelets.

All of us know that wavelet theory has made a great success in the
study of function spaces, and symbols were introduced as a
representation of operators. In this sense, operators could be
viewed as matrix under the usual wavelet bases for function
spaces, and one hopes that the above class of operators could be
characterized by the operators whose matrices under the respective
wavelet basis are privileged on the diagonal. But this is not true
except for the case where the operators themselves and their
conjugate operator all belong to Op$S^{m}_{1,1}$ (see \cite{12}).
In refs~\cite{6,17,18} one used the
Beylkin--Coifman--Meyer--Rokhlin algorithm and its generalization
to characterize the kernel-distribution of operators by their
wavelet coefficients. In analysing Calder\'on--Zygmund operators,
Yang treated their kernel-distributions as usual distribution in
$2n$ dimensions. In analysing symbol operators in Op$S^{m}_{1,
\delta}\, (0\leq \delta\leq 1)$, he treated their
kernel-distribution like distributions in $2n$ dimensions where
different coordinates play different roles. Further, one developed
pseudo-annular decomposition to study operator's continuity on the
basis of wavelet characterization (see \cite{4,11}). But there
exists difficulties to find unconditional bases for general symbol
operators in Op$S^{m}_{\rho,\delta}$ by considering their
kernel-distributions. Here, we treat directly the symbols as
distributions in phase space and {\it our first aim} is to
characterize all these symbol classes with wavelet coefficients.

Besov spaces $B^{m,\infty}_{\infty}$ is a little bigger than
H\"older spaces $C^{m}_{b}$. But the latter has no unconditional
basis, and wavelets cannot characterize it; the former has
unconditional basis, and wavelets can characterize it. Hence we
replace $S^{m}_{\rho,\delta}$ by $\tilde{S}^{m}_{\rho,\delta}$.
One writes $\sigma(x,\xi)\in \tilde{S}^{m}_{\rho,\delta}$, if
\begin{equation}
\|\partial^{\beta}_{\xi}\sigma(x,\xi)
\|_{B^{\alpha,\infty}_{\infty}} \leq C_{\alpha,\beta}
(1+|\xi|)^{m-\rho |\beta|+\delta\alpha}, \quad \forall \alpha\in
N, \beta\in N^{n}.
\end{equation}
We have the following theorem.

\begin{theor}[\!]
Given $m\in R, \rho,\delta\geq 0${\rm ,} there exists an index set
$\Lambda_{\rho,\delta}${\rm ,} a group of wavelet basis
$\{\Phi_{\lambda}(x,\xi)\} _{\lambda\in\Lambda_{\rho,\delta}}$
where $\Phi_{\lambda}(x,\xi)\in S(R^{n}\times R^{n})$ and a group
of number array spaces $N^{m}_{\rho,\delta}$ such that

\begin{enumerate}
\renewcommand\labelenumi{{\rm (\roman{enumi})}}
\leftskip .15pc
\item If $\sigma(x,\xi)\in \tilde{S}^{m}_{\rho,\delta}${\rm ,} then there
exists a unique sequence
$\{a_{\lambda}\}_{\lambda\in\Lambda_{\rho,\delta}}\in
N^{m}_{\rho,\delta}$ such that
\begin{equation*}
\hskip -1.25pc \sigma(x,\xi) = \sum\limits_{\lambda\in \Lambda_{\rho,\delta}}
a_{\lambda} \Phi_{\lambda}(x,\xi).
\end{equation*}

\item Conversely{\rm ,} if $\{a_{\lambda}\}_{\lambda\in\Lambda_{\rho,\delta}}\in
N^{m}_{\rho,\delta}${\rm ,} then there exists a unique symbol
$\sigma(x,\xi)$ such that the following formula is true in the
sense of symbol
\begin{equation*}
\hskip -1.25pc \sigma(x,\xi) = \sum\limits_{\lambda\in \Lambda_{\rho,\delta}}
a_{\lambda} \Phi_{\lambda}(x,\xi)\in \tilde{S}^{m}_{\rho,\delta}.
\end{equation*}
\end{enumerate}
\end{theor}

\begin{rem}
{\rm Given $m \in R, \rho,\delta\geq 0$, by (1.2) and (1.3), it is
easy to see that $S^{m}_{\rho,\delta} \subset
\tilde{S}^{m}_{\rho,\delta}$ where their elements are almost the
same, more precisely, $S^{m}_{\rho,\delta} {\tiny\begin{array} {c} \subset\\
\neq \end{array}} \tilde{S}^{m}_{\rho,\delta} {\tiny\begin{array}
{c} \subset\\ \neq
\end{array}} S^{m}_{\rho,\delta+\tau},\forall \tau > 0$.
Further, by the proof in Theorem 5 below, we know $S^{m}_{0,0}=
\tilde{S}^{m}_{0,0}$.}
\end{rem}

Note that, if there exist a set $S$ and a group of functions
$\{\Phi_{\lambda}(x,\xi)\}_{ \lambda\in S}$ satisfying that
$\Phi_{\lambda}(x,\xi)\in S(R^{n}\times R^{n})$ and
$\{\Phi_{\lambda}(x,\xi)\}_{\lambda\in S}$ is an orthonormal basis
in $L^{2}(R^{n}\times R^{n})$, then for each distribution
$\sigma(x,\xi)\in S'(R^{n}\times R^{n})$ and for each $\lambda\in
S$, we can define a unique number $a_{\lambda} = \langle
\sigma(x,\xi), \Phi_{\lambda}(x,\xi)\rangle$. That is to say,
there is an one-to-one relationship between the symbols in
$S'(R^{n}\times R^{n})$ and the number sequences
$\{a_{\lambda}\}_{\lambda\in S}$. Thus $\{
a_{\lambda}\}_{\lambda\in S}$ becomes a new representation for
symbol --- a wavelet representation. The difficulties to analyse
operators with wavelets are to find the appropriate wavelet basis.
The proof of Theorem~1 will be given in two sections: in Theorem~5
of \S3, we find unconditional bases for
$\tilde{S}^{m}_{0,\delta}\,(\delta\geq 0)$; in \S5, we characterize
$\tilde{S}^{m}_{\rho,\delta}\,(\rho>0)$ with wavelet coefficients.

The second generation of Calder\'on--Zygmund operators studied
kernel-distribution $k(x,z)$ where
\begin{equation*}
k(x,z) = (2\pi)^{-n}\int \sigma (x,\xi) \ \e^{iz\xi} \ \d \xi
\end{equation*}
and
\begin{equation*}
Tf(x) = \sigma (x,D)f(x) = \int k(x,z) f(x-z) \ \d z.
\end{equation*}
Meyer \cite{12} and Stein \cite{15} have established some
relations for symbol and kernel-distribution for some special
symbol class $S^{m}_{1,\delta}$. The {\it second aim} of this
paper is to get a more general result by using Theorem 1 or more
precisely, by using Theorem~6 in \S4.

\begin{theor}[\!]
Given $m\in R, \rho>0, \delta\geq 0$. If $\sigma(x,\xi)\in
S^{m}_{\rho,\delta}${\rm ,} then $\forall \alpha,\beta\in
N^{n}${\rm ,} we have
\begin{enumerate}
\renewcommand\labelenumi{{\rm (\roman{enumi})}}
\leftskip .15pc
\item If $|z|\geq \frac{1}{2}${\rm ,} then{\rm ,} $\forall \alpha, \beta\in
N^{n}${\rm ,} we have
\begin{equation*}
\hskip -1.25pc |\partial^{\alpha}_{x} \partial^{\beta}_{z} k(x,z)|\leq
C_{\alpha,\beta,N} (1+|z|)^{-N}, \forall  N > 0.
\end{equation*}

\item If $|z|\leq \frac{1}{2}${\rm ,} then{\rm ,} $\forall \alpha\in N,
\beta\in N^{n}${\rm ,} we have
\begin{equation*}
\hskip -1.25pc \|\partial^{\beta}_{z} k(x,z) \|_{B^{\alpha,\infty}_{\infty}}
\leq C_{\alpha,\beta,N} |z|^{-N}, \forall  N\geq 0
\end{equation*}
and
\begin{equation*}
\hskip -1.25pc n+m+\delta\alpha+\max(1,\rho) |\beta| < N \rho.
\end{equation*}
\end{enumerate}
\end{theor}

The proof of Theorem~{\rm 2} will be given in \S{\rm 5}.

The reason why we pay attention to the wavelet structure of
operators is to analyse precisely operator's continuity. For
example, $T1$ theorem and compensated compactness are well-known
(see \cite{2,4,9,11,16}). There are some problems which
are hard to solve without wavelets. In \cite{5,7,19}, one uses
wavelets and relative pseudo-annular decomposition to study the
$T1$ theorem and the compensated theory and gets some good
results. The {\it third aim} of this paper is to study the
$L^{2}$-continuity and the $L^{p}$-continuity of non-regular symbol
operators. On the basis of symbol's wavelet coefficients in phase
space, we can apply a precise Huygens' principal (or a precise
micro-analysis method) to study operator's continuity (see also \cite{10}).

In \cite{15}, Stein studied the $L^{2}$-continuity of operators
defined by the symbol in $S^{0}_{0,0}= C^{\infty}_{b}(R^{2n})$. In
\cite{3} and \cite{8}, one studied pseudodifferential operators in
phase space. In \cite{1}, Calder\'on and Vaillancourt studied
$L^{2}$-continuity of symbol operators where symbol
$\sigma(x,\xi)$ belong to the H\"older space
$C^{2n+1}_{b}(R^{2n})$ and in some sense, which is the special
Besov space $ B^{2n+1,\infty}_{\infty}(R^{2n}) =
B^{2n+1,\infty}_{\infty}$. Here we reduce an index $n$ for the
order of smoothness and establish $L^{2}$-continuity also; in
fact, for $s\leq n<s'$, we know that
$B^{s',\infty}_{\infty}\subset B^{n,1}_{\infty}\subset
B^{s,\infty}_{\infty}$. Further, we can construct a special
operator to show that our result is sharp. That is our Theorem~3.

\begin{theor}[\!]$\left.\right.$

\begin{enumerate}
\renewcommand\labelenumi{{\rm (\roman{enumi})}}
\leftskip .15pc
\item If $\sigma(x,\xi)\in B^{n,1}_{\infty}${\rm ,} then we have
\begin{equation}
\hskip -1.25pc \sigma(x,D) \ \hbox{defines an operator which is continuous from}
\ L^{2}(R^{n}) \ \mbox{to} \ L^{2}(R^{n}).
\end{equation}

\item Conversely{\rm ,} for $0 < s < n${\rm ,} there exists a symbol
$\sigma(x,\xi)\in B^{s,\infty}_{\infty}$ but
\begin{equation}
\hskip -1.25pc \sigma(x,D) \ \hbox{is not continuous from} \ L^{2}(R^{n}) \
\hbox{to} \ L^{2}(R^{n}).
\end{equation}
\end{enumerate}
\end{theor}

In addition, if we strengthen a little the above assumption, we
can consider $L^{p}$-continuity. Let $Q = \{x = (x_{1},\ldots,
x_{n}), 0\leq x_{i}\leq 1, 1 \leq i\leq n\}$ be a unit cube in
$R^{n}$. For $j\geq 0$ and $k\in Z^{n}$, denote
$2^{-j}k+2^{-j}Q=\{x\hbox{:}\ 2^{j}x-k\in Q\}$. Let $I_{n}$ be the set
which is composed by $n$ elements in $R^{n}$ which are the unit
vectors in the direction of the axes. For arbitrary distribution
$f(x)$ and for $e\in I_{n}, h\in R, m\in N$, let
\begin{equation*}
\tau_{he}f(x) = f(x+he)-f(x) \quad \hbox{and} \quad \tau^{m}_{he}
= (\tau_{he})^{m}.
\end{equation*}
For $j\geq 1, X=(x,\xi)\in R^{2n}, \ \e \in I_{2n}$, denote
\begin{equation*}
\sigma_{j,e}(X) = \tau^{n}_{2^{-j}e} \sigma(X).
\end{equation*}
Denote
\begin{equation*}
\omega(0) = \sup\limits_{k\in
Z^{n}}\int_{k+Q}\!\int_{R^{n}}|\sigma(x,\xi)|\ \d x \ \d \xi,
\end{equation*}
and for $j\geq 1$, denote
\begin{equation*}
\omega(j)=\sup\limits_{k\in Z^{n}, e\in I_{2n}}
\int_{2^{-j}k+2^{-j}Q}\! \int_{R^{n}} |\sigma_{j,e}(x,\xi)| \ \d x
\ \d \xi.
\end{equation*}
We say that $\sigma(x,\xi) \in B^{s}$, if
\begin{equation*}
\sum\limits_{j} 2^{(n+s)j} \omega(j) < \infty.
\end{equation*}
By (2.7) and (7.1) below, we know that $B^{n}{\tiny \begin{array}{c}\subset\\
\neq\end{array}} B^{n,1}_{\infty}$. Now we establish
$L^{p}$-continuity.

\begin{theor}[\!]$\left.\right.$

\begin{enumerate}
\renewcommand\labelenumi{{\rm (\roman{enumi})}}
\leftskip .15pc
\item If $\sigma(x,\xi)$ satisfies the condition
\begin{equation}
\hskip -1.25pc \sigma(x,\xi) \in B^{n},
\end{equation}
then for $1\leq p\leq\infty${\rm ,} we have
\begin{equation}
\hskip -1.25pc \sigma(x,D) \ \hbox{is continuous from} \
L^{p}(R^{n}) \ \mbox{to} \  L^{p}(R^{n}).
\end{equation}

\item Conversely{\rm ,} for $0 < s < n${\rm ,} there exists $\sigma(x,\xi)$
satisfies the condition
\begin{equation}
\hskip -1.25pc \sigma(x,\xi) \in B^{s},
\end{equation}
but for $1\leq p\leq\infty${\rm ,} we have
\begin{equation}
\hskip -1.25pc \sigma(x,D) \ \hbox{is not continuous from} \ L^{p}(R^{n}) \
\hbox{to} \ L^{p}(R^{n}).
\end{equation}
\end{enumerate}
\end{theor}

The difficulty to study operator's continuity is to find an
appropriate operator's decomposition such that the relative
operators have some pseudo-orthogonality. Our new idea is to
establish the operator's continuity in Theorems~3 and 4 on the
basis of wavelet characterization in phase space, and the proof
will be given in the last two sections of this paper.

\begin{rem}

{\rm Meyer wrote in his famous book \cite{12} that, for a long
time, the study of operators stayed in two isolated
classes---Calder\'on--Zygmund operators and symbol operators. On one
hand, one has found wavelet characterization for
Calder\'on--Zygmund operators and established such operator's
continuity and also commutator operator's continuity (see
\cite{5,7,13,17,19}). On the other hand, one has given
a wavelet representation for symbol operators and developed
relative methods to study operator's continuity in this paper and
in other papers (see \cite{18}). That is to say, we can study both
Calder\'on--Zygmund operators and symbol operators under their
wavelet representation.}
\end{rem}

\section{Preliminaries}
\setcounter{equation}{0}

At the begin of this section, we introduce some notations for
wavelets and prove some wavelet properties.

In this paper, we use a wavelet basis which is a tensor product of
the wavelets in dimension 1. When we characterize
$S^{m}_{\rho,\delta}$ in \S\S3 and 4 and when we analyse
kernel-distribution in \S5, we always use Meyer's wavelets; but
in other cases, when we prove operator's continuity, we need
Meyer's wavelets; when we construct special operators to prove
that our results are sharp, we need sufficiently regular
Daubechies' wavelets. In dimension 1, denote the father wavelet by
$\Phi^{0}(x)$ and the mother wavelet by $\Phi^{1}(x)$. In high
dimension, for $\epsilon= (\epsilon_{1},\ldots,\epsilon_{n}) \in
\{0,1\}^{n}$, denote
\begin{equation}
\Phi^{\epsilon}(x)=\prod\limits^{n}_{i=1}
\Phi^{\epsilon_{i}}(x_{i}) \quad \hbox{and} \quad
\Phi^{0}(x)=\Phi^{(0,\ldots, 0)}(x).
\end{equation}
For $j\in Z, k\in Z^{n}$, denote
\begin{equation}
f_{j,k}(x)=2^{nj/2}f(2^{j}x-k).
\end{equation}
Let $\{V_{j}\}_{j\in Z}$\vspace{.05pc} be an orthogonal multi-resolution
analysis in $L^{2}(R^{n})$ and $V_{j+1}=V_{j}\oplus W_{j}$. Then
$\{\Phi^{0}_{j,k}(x)\}_{k\in Z^{n}}$ is an orthonormal wavelet
basis in $V_{j}$ and
$\{\Phi^{\epsilon}_{j,k}(x)\}_{\epsilon\in\{0,1\}^{n}\backslash
\{0\}, k\in Z^{n}}$ is an orthonormal wavelet basis in $W_{j}$ and
$L^{2}(R^{n})=V_{0}\bigoplus_{j\geq 0}W_{j}$. Let $P_{j}$
be the projector operator from $L^{2}(R^{n})$ to $V_{j}$ and let
$Q_{j}$ be the projector operator from $L^{2}(R^{n})$ to $W_{j}$.
It is easy to see that $P_{0}+\sum_{j\geq 0} Q_{j}$ is the
unit operator $I$. Let
\begin{align}
\Lambda_{n} &= \{ \lambda = (\epsilon, j, k), \epsilon\in
\{0,1\}^{n}, j\geq 0, k\in Z^{n};\nonumber\\[.3pc]
&\quad\,\hbox{and if} \ j>0, \ \hbox{then} \ \epsilon\neq 0\}.
\end{align}
Then $\{\Phi^{\epsilon}_{j,k}(x)\}_{(\epsilon,j,k)\in\Lambda_{n}}$\vspace{.05pc}
is an orthonormal wavelet basis in $L^{2}(R^{n})$. $\forall
\epsilon\in \{0,1\}^{n}$, there exists $\{g^{\epsilon}_{k}\}_{k\in
Z^{n}}$ such that
\begin{equation}
\Phi^{\epsilon}(x) = \sum\limits_{k} g^{\epsilon}_{k}
\Phi^{0}(2x-k).
\end{equation}
$\forall \epsilon = (\epsilon_{1},\ldots,\epsilon_{n})\neq 0$, let
$\tau_{\epsilon}$ denote the smallest number $i$ such that
$\epsilon_{i}\neq 0$ and let $e_{\epsilon} = e_{\tau_{\epsilon}}$
denote the vector where the $\tau_{\epsilon}$ coordinate is 1 and
the rest are 0. For any sequence $\{a_{k}\}_{k\in Z_{n}}$, %%abi 6-10
let $\tau^{0}=S^{0}$ be the unit operator satisfying $\tau^{0} a_{k} =
S^{0}a_{k} =a_{k}$; for $ e_{i}\in I_{n}$ where its $i$th element
is 1 and the rest are $0$; and for $s\in N$, let
\begin{equation}
\tau_{\pm e_{i}} a_{k} = a_{k\pm e_{i}}-a_{k}\quad \hbox{and}\quad
S_{e_{i}} a_{k} = - \sum\limits ^{-1+k_{i}} _{l=-\infty}
a_{(k_{1},\ldots, k_{i-1},l,k_{i+1},\ldots, k_{n})},
\end{equation}
and let $\tau_{\pm e_{i}}^{s}=(\tau_{\pm e_{i}})^{s}$ and
$S_{e_{i}}^{s} = (S_{e_{i}})^{s}$. Further, for $\alpha\in N^{n}$,
let
\begin{equation}
\tau^{\alpha}_{\pm}=\prod\limits^{n}_{i=1}
\tau^{\alpha_{i}}_{\pm e_{i}} \quad\hbox{and}\quad S^{\alpha} =
\prod\limits^{n}_{i=1} S^{\alpha_{i}}_{ e_{i}}.
\end{equation}
For $e_{i}\in I_{n}$ such that the $i$th element of $e_{i}$ is 1
and for $s\in N$, let $S_{e_{i}}^{0} f(x)= f(x)$ and $S_{e_{i}}
f(x) = - \sum_{l=-\infty}^{-1} f(x-l e_{i})$, and let
$S_{e_{i}}^{s} = (S_{e_{i}})^{s}$; further, for $\alpha\in N^{n}$,
let $S^{\alpha} = \prod_{i=1}^{n} S^{\alpha_{i}}_{e_{i}}$.

\begin{lem}$\left.\right.$
\begin{enumerate}
\renewcommand\labelenumi{{\rm (\roman{enumi})}}
\leftskip .1pc
\item For $\epsilon\in \{0,1\}^{n}\backslash 0$ and $s\in N${\rm
,} $\tilde{\Phi}^{\epsilon,s}(x)= \sum_{k}
(S^{s}_{e_{\epsilon}}g^{\epsilon}_{k})\Phi^{0}(2x-k)$ satisfies
$\Phi^{\epsilon}(x) = \tau_{-\frac{1}{2}e_{\epsilon}}^{s}
\tilde{\Phi}^{\epsilon,s}(x)${\rm ;} and further{\rm ,} if
$\Phi^{\epsilon}(x)$ are Meyer{\rm '}s wavelets{\rm ,} then
$\tilde{\Phi}^{\epsilon,s}(x)\in S(R^{n})${\rm ;} if
$\Phi^{\epsilon}(x)$ are Daubechies{\rm '} wavelets and $s$ is less than
the index of divergence moment of wavelets{\rm ,} then
$\tilde{\Phi}^{\epsilon,s}(x)$ have compact support.

\item For Meyer{\rm '}s wavelet{\rm ,} $\forall \beta\in N^{n}${\rm ,} $S^{\beta}
(\partial^{\beta}\Phi^{0})(x)\in S(R^{n})$. \end{enumerate}
\end{lem}

\begin{proof}$\left.\right.$
\begin{enumerate}
\renewcommand\labelenumi{(\roman{enumi})}
\leftskip .1pc
\item For $\epsilon\in \{0,1\}^{n}\backslash 0$ and $s\in N$, by the
scale equation $\Phi^{\epsilon}(x) = \sum_{k} g^{\epsilon}_{k}
\Phi^{0}(2x-k)$ and by the construction of
$\tilde{\Phi}^{\epsilon,s}(x)$, we have $\Phi^{\epsilon}(x) =
\tau_{-\frac{1}{2}e_{\epsilon}}^{s} \tilde{\Phi}^{\epsilon,s}(x)$.
Further, by divergence moment properties of wavelets, we have:

\bigskip
%\begin{enumerate}
%\renewcommand\labelenumii{(\arabic{enumii})}
(1)\ \ If $\Phi^{\epsilon}(x)$ are Meyer's wavelets, then
$|S^{s}_{e_{\epsilon}}g^{\epsilon}_{k}|\leq C_{s,N}
(1+|k|)^{-N},\forall N>0$ and hence
$\tilde{\Phi}^{\epsilon,s}(x)\in S(R^{n})$.

(2)\ \ If $\Phi^{\epsilon}(x)$ are Daubechies' wavelets and $s$ is
less than the index of divergence moment of wavelets, then there
exists $C_{s}$ such that, for $|k|\geq C_{s}$,
$S^{s}_{e_{\epsilon}}g^{\epsilon}_{k}=0$, and hence
$\tilde{\Phi}^{\epsilon,s}(x)$ have compact support.\vspace{-.5pc}
%\end{enumerate}
\bigskip

\item For Meyer's wavelet, $\forall \beta\in N^{n}$, $\sum_{k}
k^{\beta} \Phi^{0}(x-k)$ are polynomials $P_{\beta}(x)$ where the
degree of $x_{i}$ is $\beta_{i}$; hence we have $S^{\beta}
(\partial^{\beta}\Phi^{0})(x)\in S(R^{n})$. Or we can prove (ii)
by the fact that the Fourier transformation of $S^{\beta}
(\partial^{\beta}\Phi^{0})(x)$ is equal to
$f_{\beta}(\xi)=C_{\beta}\xi^{\beta}\prod^{n}_{j=1} (1-
\hbox{e}^{i \xi_{j}} )^{-\beta_{j}} \hat{\Phi}^{0}(\xi)$; since
${\rm supp } \hat{\Phi}^{0}(\xi)\subset [-\frac{4\pi}{3},
\frac{4\pi}{3}]^{n}$, so $f_{\beta}(\xi) \in
S(R^{n})$.\vspace{-1.5pc}
\end{enumerate}
\end{proof}

Besov spaces $B^{s,q}_{p}$, which were introduced systematically
by Peetre \cite{14} can be characterized with their wavelet
coefficients (see \cite{12} and \cite{19}). For $f(x) =
\sum_{\lambda=(\epsilon,j,k)\in\Lambda_{n}}
a_{\lambda}\Phi_{\lambda}(x)$, we have the following.

\begin{lem}$\left.\right.$\vspace{-1.7pc}

\begin{equation}
f(x)\in B^{s,q}_{p}(R^{n})\Longleftrightarrow
\left(\sum\limits_{j\geq 0} 2^{jq(s+\frac{n}{2}-\frac{n}{p})}
\left(\sum\limits_{\epsilon,k}
|a_{\lambda}|^{p}\right)^{\frac{q}{p}}\right)^{\frac{1}{q}}
<\infty.
\end{equation}
\end{lem}

Secondly, we introduce a useful and simple inequality which would
be used often in this paper.

\begin{lem}
$\forall \alpha\geq 1, m\in R, x,y\in R^{n}${\rm ,} we have
\begin{equation}
(1+|x|)^{m}\leq (1+|y|)^{m} (1+\alpha |x-y|)^{|m|}.
\end{equation}
\end{lem}

\begin{proof}\vspace{.2pc}
It is evident for $m\geq 0$. If $m<0$, then we have
\begin{align*}
(1+|x|)^{m}\leq (1+|y|)^{m} (1+|x-y|)^{|m|}\leq (1+|y|)^{m}
(1+\alpha |x-y|)^{|m|}.
\end{align*}

At the end of this section, a variation of the result in \cite{15}
(which discusses the operator's continuity) will be introduced.
For $j\geq 0$ and $m=(k,l)\in Z^{2n}$, let $T^{*}_{j,m}$ be the
relative conjugate operators of operators $T_{j,m}$. Then we have
the following lemma.
\end{proof}

\begin{lem}
Suppose that $T_{j,m}$ satisfies the following three
conditions{\rm :}
\begin{align}
\|T_{j,m} \|_{L^{2}\rightarrow L^{2}} &\leq C,\\[.3pc]
\|T_{j,k_{1},l_{1}}T^{*}_{j,k_{2},l_{2}} \|_{L^{2}\rightarrow
L^{2}} &\leq C(1+4^{-j}|k_{1}-k_{2}|)^{-2N_{0}}
(1+|l_{1}-l_{2}|)^{-2N_{0}},\\[.3pc]
\|T^{*}_{j,k_{1},l_{1}}T_{j,k_{2},l_{2}} \|_{L^{2}\rightarrow
L^{2}} &\leq C(1+|k_{1}-k_{2}|)^{-2N_{0}}
(1+4^{-j}|l_{1}-l_{2}|)^{-2N_{0}}.
\end{align}
Then for $N_{0}>n${\rm ,} $T_{j}=\sum_{m\in Z^{2n}} T_{j,m}$
defines an operator which is continuous from $L^{2}$ to $L^{2}$
and $\|T_{j}\|_{L^{2}\rightarrow L^{2}}\leq C4^{jn}$.
\end{lem}

\begin{proof}
First, we consider a finite sum $S_{j}=S_{j,N}=\sum_{|m|\leq N}
T_{j,m}$. Since $S_{j}^{*}S_{j}$ is a self-adjoint operator, we
have
$\|S_{j}\|^{2}=\|S_{j}^{*}S_{j}\|=\|(S_{j}^{*}S_{j})^{M}\|^{1/M}$
for all integer $M$. But we have\vspace{-.2pc}
\begin{align}
(S_{j}^{*}S_{j})^{M}&=\sum\limits_{k_{1},l_{1}}
\sum\limits_{k_{2},l_{2}}\cdots \sum\limits_{k_{2M-1},l_{2M-1}}\nonumber\\[.3pc]
&\quad\,\times\sum\limits_{k_{2M},l_{2M}}
T^{*}_{j,k_{1},l_{1}}T_{j,k_{2},l_{2}}\cdots
T^{*}_{j,k_{2M-1},l_{2M-1}}T_{j,k_{2M},l_{2M}}.
\end{align}
We maximize $\|(S_{j}^{*}S_{j})^{M}\|$ by
\begin{equation}
\sum\limits_{k_{1},l_{1}}
\sum\limits_{k_{2},l_{2}}\cdots \sum\limits_{k_{2M-1},l_{2M-1}}
\sum\limits_{k_{2M},l_{2M}}
\|T^{*}_{j,k_{1},l_{1}}T_{j,k_{2},l_{2}}\cdots
T^{*}_{j,k_{2M-1},l_{2M-1}}T_{j,k_{2M},l_{2M}}\|.
\end{equation}

First, we re-group all the operators two by two, and apply the
continuity of $\|T^{*}_{j,m}T_{j,m'}\|$. We get\vspace{-.2pc}
\begin{align*}
&\|T^{*}_{j,k_{1},l_{1}}T_{j,k_{2},l_{2}}\cdots
T^{*}_{j,k_{2M-1},l_{2M-1}}T_{j,k_{2M},l_{2M}}\|\\[.3pc]
&\quad\,\leq C^{M}
(1+|k_{1}-k_{2}|)^{-2N_{0}}(1+4^{-j}|l_{1}-l_{2}|)^{-2N_{0}}\\[.3pc]
&\qquad\,\times \cdots \times (1+|k_{2M-1}-k_{2M}|)^{-2N_{0}}
(1+4^{-j}|l_{2M-1}-l_{2M}|)^{-2N_{0}}.
\end{align*}
\pagebreak

\noindent Then we maximize $\|T^{*}_{j,k_{1},l_{1}}\|$ and
$\|T_{j,k_{2M},l_{2M}}\|$ by the constant $C$, then re-group the
remaining operators two by two. Applying the continuity of
$\|T_{j,m}T^{*}_{j,m'}\|$, we get
\begin{align*}
&\|T^{*}_{j,k_{1},l_{1}}T_{j,k_{2},l_{2}}\cdots
T^{*}_{j,k_{2M-1},l_{2M-1}}T_{j,k_{2M},l_{2M}}\| \\[.3pc]
&\quad\,\leq C^{M+2}
(1+4^{-j}|k_{2}-k_{3}|)^{-2N_{0}}(1+|l_{2}-l_{3}|)^{-2N_{0}}\\[.3pc]
&\qquad\,\times \cdots \times (1+4^{-j}|k_{2M-2}-k_{2M-1}|)^{-2N_{0}}
(1+|l_{2M-2}-l_{2M-1}|)^{-2N_{0}}.
\end{align*}
Combining the above two cases, we have
\begin{align*}
&\|T^{*}_{j,k_{1},l_{1}}T_{j,k_{2},l_{2}}\cdots
T^{*}_{j,k_{2M-1},l_{2M-1}}T_{j,k_{2M},l_{2M}}\| \leq C^{M+1}
(1+|k_{1}-k_{2}|)^{-N_{0}}\\[.3pc]
&\quad\,\times (1+4^{-j}|k_{2}-k_{3}|)^{-N_{0}} \cdots
(1+|k_{2M-1}-k_{2M}|)^{-N_{0}}\\[.3pc]
&\quad\,\times
(1+4^{-j}|l_{1}-l_{2}|)^{-N_{0}}(1+|l_{2}-l_{3}|)^{-N_{0}}\\[.3pc]
&\quad\,\times \cdots \times (1+4^{-j}|l_{2M-1}-l_{2M}|)^{-N_{0}}.
\end{align*}
Summing in order $k_{1},\ldots,k_{2M-1}$ and $l_{1},\ldots,
l_{2M-1}$, one gets $C^{M+1} 4^{jn(2M-1)}$; then summing $k_{2M}$
and $l_{2M}$, one gets
\begin{equation}
\|S\|^{2M}\leq CN^{2n}C^{M+1} 4^{jn(2M-1)} \quad\hbox{or}\quad
\|S\|\leq (CN^{2n}C^{M+1} 4^{jn(2M-1)})^{\frac{1}{2M}}.
\end{equation}
Letting $M\rightarrow\infty$, we get $\|S\|\leq C4^{jn}$.

Further, we adopt Journ\'e's methods to pass to the general case.
According to the above result, $\forall f(x)\in L^{2}$, we have
\begin{equation}
\left\|\sum\limits_{|m|\leq N} \lambda_{m}T_{m}f(x)\right\|_{L^{2}}
\leq C\|f(x)\|_{L^{2}},\quad \forall N\in N, |\lambda_{m}|\leq 1.
\end{equation}
Let $\epsilon(N)=\sup_{\tilde{N}\geq N} \|\sum_{N\leq |m|\leq
\tilde{N}} \lambda_{m} T_{m}f(x)\|_{L^{2}}$. To prove that
$\sum_{m} T_{m}f(x)$ converges to a function in $L^{2}$, it is
sufficient to prove that $\lim_{N\rightarrow\infty}\epsilon(N)=0$.
It is evident that, $\forall N\leq N'$, we have $\epsilon(N)\geq
\epsilon(N')$. If $\epsilon(N)$ does not approach zero, then there
exists $\delta>0$ and $N>0$ such that $\epsilon (N') \geq \delta,
\forall N'\geq N$. Then we can choose
$m^{1}_{N}<m^{2}_{N}<\cdots<m^{2k}_{N}<m^{2k+1}_{N}<\cdots$ such
that for $Z_{k}=\sum_{ m^{2k}_{N}\leq |m| \leq m^{2k+1}_{N}}
T_{m}f(x)$, we have
\begin{equation}
\|Z_{k}\|_{L^{2}}\geq \delta.
\end{equation}

For $\theta=(\theta_{1},\ldots,\theta_{k})\in\{-1,1\}^{k}$, let
$Z(\theta,k)= \sum^{k}_{i=1} \theta_{i} Z_{i}$. According to
(2.15), we have $\| Z(\theta,k)\|\leq C\|f(x)\|_{L^{2}}$. Since
$\sum^{k}_{i=1}\|Z_{i}\|^{2}_{L^{2}}\leq 2^{-k}\sum_{\theta\in
\{-1,1\}^{k}} \|Z(\theta,k)\|^{2}$, we have
$\sum^{k}_{i=1}\|Z_{i}\|_{L^{2}}\leq C\|f\|_{L^{2}}$, which
contradicts (2.16)!
\end{proof}

\section{Unconditional bases for $\pmb{\tilde{S}^{m}_{0,\delta}\,(\delta\geq 0)}$}

In this section, we use the usual $2n$ dimension wavelet basis in
phase space to characterize $S^{m}_{0,0}=\tilde{S}^{m}_{0,0}$ and
use wavelet basis which comes from tensor product of wavelet basis
in $n$ dimension to characterize
$\tilde{S}^{m}_{0,\delta}\,(\delta>0)$ .

Let $\Lambda_{0,0}=\Lambda_{2n}$ and $\forall \lambda=
(\epsilon,\epsilon',j,k,l)\in \Lambda_{0,0}$, let
$\Phi_{\lambda}(x,\xi) = \Phi^{\epsilon}_{j,k}(x)
\Phi^{\epsilon'}_{j,l}(\xi)$. Then
$\{\Phi_{\lambda}(x,\xi)\}_{\lambda\in\Lambda_{0,0}}$ is an
orthonormal basis in $L^{2}(R^{n}\times R^{n})$. For $a_{\lambda}
= \langle \sigma(x,\xi), \Phi_{\lambda}(x,\xi)\rangle $, the
following equality is true in the sense of distribution:
\setcounter{equation}{0}
\begin{equation}
\sigma(x,\xi)=\sum\limits_{\lambda\in\Lambda_{0,0}} a_{\lambda}
\Phi_{\lambda}(x,\xi).
\end{equation}
Hence, we know that $\{ a_{\lambda} \}_{\lambda\in \Lambda_{0,0}}$
becomes a new representation for symbol $\sigma(x,\xi)$. We say
that $\{ a_{\lambda} \}_{\lambda\in \Lambda_{0,0}}\in
N^{m}_{0,0}$, if
\begin{equation}
|a_{\lambda}|\leq C_{N} 2^{-jN} (1+|2^{-j}l|)^{m},\quad \forall N>0,
\lambda\in \Lambda_{0,0}.
\end{equation}
For $\delta>0$, let $\Lambda_{0,\delta}= \Lambda_{n}\times
\Lambda_{n}$; and for $\lambda=(\epsilon,j,k,\epsilon',j',k')\in
\Lambda_{0,\delta}$, let $\Phi_{\lambda}(x,\xi) =
\Phi^{\epsilon}_{j,k}(x) \Phi^{\epsilon'}_{j',k'}(\xi)$. Then
$\{\Phi_{\lambda}(x,\xi)\}_{\lambda\in \Lambda_{0,\delta}}$ is an
orthonormal wavelet basis in $L^{2}(R^{n}\times R^{n})$. For
$a_{\lambda}=\langle \sigma(x,\xi),\Phi_{\lambda}(x,\xi)\rangle$,
it is clear that $\{ a_{\lambda}\}_{\lambda\in
\Lambda_{0,\delta}}$ becomes a new representation for symbol. For
$\delta>0$, we write $\{ a_{\lambda}\}_{\lambda\in
\Lambda_{0,\delta}}\in N^{m}_{0,\delta}$, if
\begin{equation}
|a_{\lambda}| \leq
 C_{\alpha,\beta} 2^{-(\frac{n}{2}+\alpha)j}  2^{-(\frac{n}{2}+\beta)j'}
(1+|2^{-j'}k'|) ^{m+\delta\alpha},\quad \forall \alpha,\beta\geq 0.
\end{equation}
On basis of the above notation, for $\delta\geq 0$, we have the
following.

\begin{theor}[\!]
The following two conditions are equivalent{\rm :}
\begin{align}
\sigma(x,\xi) &\in \tilde{S}^{m}_{0,\delta},\\[.3pc]
\{a_{\lambda}\}_{\lambda \in \Lambda_{0,\delta}} &\in
N^{m}_{0,\delta}.
\end{align}
\end{theor}

\begin{proof}
{\it First step.} We consider the case where $\delta=0$ and we
prove that $\sigma(x,\xi)\in \tilde{S}^{m}_{0,0}$ implies that $\{
a_{\lambda} \}_{\lambda\in \Lambda_{0,0}}\in N^{m}_{0,0}$. We
consider three cases: (i) $\epsilon'\neq 0$, (ii) $\epsilon'= 0,
\epsilon\neq 0$ and (iii) $\epsilon=\epsilon'= 0$. For arbitrary
$\epsilon\in \{0,1\}^{n}\backslash \{0\}$ and $N>0$, let
$I^{N}_{\epsilon}f(x)$ be the $N$th integration of $f(x)$ for the
$\tau_{\epsilon}$-coordinate. For Case (i) and for sufficiently
large $N'>2N+n+|m|$, we have
\begin{align*}
|a_{\lambda}| &= |\langle \sigma (x,\xi), \Phi^{\epsilon}_{j,k}(x)
\Phi^{\epsilon'}_{j,l}(\xi)\rangle |\\[.3pc]
&= 2^{j(n-N)} |\langle
\partial^{N}_{\xi_{\tau_{\epsilon'}}}\sigma (x,\xi),
\Phi^{\epsilon}(2^{j}x-k)
(I^{N}_{\epsilon'}\Phi)^{\epsilon'}(2^{j}\xi-l)\rangle|\\[.3pc]
&\leq 2^{j(n-N)} \int |\langle
\partial^{N}_{\xi_{\tau_{\epsilon'}}}\sigma (x,\xi),
\Phi^{\epsilon}(2^{j}x-k)\rangle |
|(I^{N}_{\epsilon'}\Phi)^{\epsilon'}(2^{j}\xi-l) | \d\xi\\[.3pc]
&\leq C 2^{-jN} \int \frac {(1+|\xi|)^{m}} {(1+|2^{j}\xi-l|)^{N'}}
\d\xi.
\end{align*}
Then applying Lemma 3 to $(1+|\xi|)^{m}$, we have
\begin{align*}
|a_{\lambda}| &\leq C 2^{-jN} (1+|2^{-j}l|)^{m} \int
(1+ |2^{j}\xi-l|) ^{|m|-N'} \d\xi \\[.3pc]
&\leq  C 2^{-j(n+N)} (1+|2^{-j}l|)^{m}.
\end{align*}
For Case (ii), by Lemma 2, we have
\begin{align*}
|a_{\lambda}| =  \int |\langle \sigma (x,\xi),
\Phi^{\epsilon}_{j,k}(x)\rangle | |\Phi^{\epsilon'}_{j,l}(\xi) |
\d\xi \leq   C 2^{-jN} \int   \frac {(1+|\xi|)^{m}}
{(1+|2^{j}\xi-l|)^{N'}}   \d\xi.
\end{align*}
Then applying Lemma 3 to $(1+|\xi|)^{m}$, we have
\begin{align*}
|a_{\lambda}| &\leq C 2^{-jN} (1+|2^{-j}l|)^{m} \int (1+
|2^{j}\xi-l|)^{|m|-N'} \d\xi\\[.3pc]
&\leq C 2^{-j(n+N)} (1+|2^{-j}l|)^{m}.
\end{align*}
For Case (iii), by Lemmas 2 and 3, we have
\begin{align*}
|a_{\lambda}| &= \int |\langle \sigma (x,\xi),
\Phi^{0}(x-k)\rangle | |\Phi^{0}(\xi-l) |
\d\xi \leq C \int \frac {(1+|\xi|)^{m}} {(1+|\xi-l|)^{N'}} \d\xi\\[.3pc]
&\leq C (1+|l|)^{m} \int (1+ |\xi-l|) ^{|m|-N'} \d\xi \leq C
(1+|l|)^{m}.
\end{align*}
{\it Second step.} We consider the case where $\delta=0$ and we
prove that $\{ a_{\lambda} \}_{\lambda\in \Lambda_{0,0}}\in
N^{m}_{0,0}$ implies that $\sigma(x,\xi)\in S^{m}_{0,0}$. For
arbitrary $\alpha,\beta\in N^{n}$, we choose $N>
n+|\alpha|+|\beta|$ and $N'>n+|m|$. We have
\begin{align*}
\Bigg|\partial^{\alpha}_{x} \partial^{\beta}_{\xi}
\sum\limits_{\lambda\in \Lambda_{0,0}} a_{\lambda}
\Phi_{\lambda}(x,\xi) \Bigg| &\leq \sum\limits_{j\geq 0}
2^{j(n+|\alpha|+|\beta|-N)} \sum\limits_{k} |(\partial
^{\alpha}_{x} \Phi^{\epsilon})(2^{j}x-k)|\\[.3pc]
&\quad\, \times\sum\limits_{l} (1+|2^{-j}l|)^{m}
|(\partial ^{\beta}_{\xi} \Phi^{\epsilon'})(2^{j}\xi-l)| \\[.3pc]
&\leq \sum\limits_{j\geq 0} 2^{j(n+|\alpha|+|\beta|-N)}
\sum\limits_{l} \frac{(1+|2^{-j}l|)^{m}} {(1+|2^{j}\xi-l|)^{-N'}}.
\end{align*}
Then applying Lemma 3 to $(1+|2^{-j}l|)^{m}$, we have
\begin{align*}
\Bigg|\partial^{\alpha}_{x} \partial ^{\beta}_{\xi}
\sum\limits_{\lambda\in \Lambda_{0,0}} a_{\lambda}
\Phi_{\lambda}(x,\xi)\Bigg| &\leq C_{\alpha,\beta} (1+|\xi|)^{m}
\sum\limits_{j\geq 0} 2^{j(n+|\alpha|+|\beta|-N)}\\[.3pc]
&\leq C_{\alpha,\beta} (1+|\xi|)^{m}.
\end{align*}
{\it Third step.} We consider the case where $\delta>0$ and we
prove that (3.4) implies (3.5). We distinguish four cases.
\begin{enumerate}
\renewcommand\labelenumi{(\arabic{enumi})}
\leftskip .1pc
\item If $\epsilon=\epsilon'=0$, then
$a^{\epsilon,\epsilon'}_{j,k,j',k'}= a^{0,0}_{0,k,0,k'} = \langle
\sigma (x,\xi), \Phi^{0}(x-k)\Phi^{0}(\xi-k')\rangle $. Since
$|\langle \sigma(x,\xi), \Phi^{0}(x-k)\rangle|\leq
C(1+|\xi|)^{m}$, we apply Lemma 3, and get
\begin{equation*}
\hskip -1.25pc |a^{0,0}_{0,k,0,k'}|\leq C (1+|k'|)^{m}.
\end{equation*}

\item If $\epsilon=0,\epsilon'\neq 0$, then
\begin{align*}
\hskip -1.25pc a^{\epsilon,\epsilon'}_{j,k,j',k'} & =
a^{0,\epsilon'}_{0,k,j',k'}= \langle \sigma (x,\xi), \Phi^{0}(x-k)
\Phi^{\epsilon}_{j',k'}(\xi)\rangle\\[.3pc]
&= 2^{-j'|\beta|} \langle
\partial^{\beta}_{\xi_{\tau_{\epsilon}}} \sigma (x,\xi),
\Phi^{0}(x-k) (I^{\beta}_{\epsilon}
\Phi^{\epsilon})_{j',k'}(\xi)\rangle .
\end{align*}
Since
$|\langle\partial^{\beta}_{\xi_{\tau_{\epsilon}}}\sigma(x,\xi),
\Phi^{0}(x-k)\rangle |\leq C(1+|\xi|)^{m}$, we apply Lemma~3, and
get
\begin{equation*}
\hskip -1.25pc |a_{\lambda}| \leq C_{N} 2^{-j'N} (1 + |2^{-j'}k' |)^{m}.
\end{equation*}

\item If $\epsilon\neq 0, \epsilon'=0$, then
\begin{equation*}
\hskip -1.25pc a^{\epsilon,\epsilon'}_{j,k,j',k'}= a^{\epsilon,0}_{j,k,0,k'} =
\langle \sigma (x,\xi), \Phi^{\epsilon}_{j,k}(x)
\Phi^{0}(\xi-k')\rangle .
\end{equation*}
Hence by Lemmas 2 and 3, we get
\begin{equation*}
\hskip -1.25pc |a^{\epsilon,0}_{j,k,0,k'}|\leq C
2^{-(\frac{n}{2}+\alpha)j}(1+|k'|)^{m+\delta\alpha}.
\end{equation*}

\item If $|\epsilon||\epsilon'|\neq 0$, then
\begin{align*}
\hskip -1.25pc a^{\epsilon,\epsilon'}_{j,k,j',k'} &= \langle \sigma (x,\xi),
\Phi^{\epsilon}_{j,k}(x)\Phi^{\epsilon'}_{j',k'}(\xi)\rangle\\[.3pc]
&= 2^{-j'|\beta|} \langle
\partial^{\beta}_{\xi_{\tau_{\epsilon}}} \sigma (x,\xi),
\Phi^{\epsilon}_{j,k}(x) (I^{\beta}_{\epsilon'}
\Phi^{\epsilon'})_{j',k'}(\xi)\rangle .
\end{align*}
Hence by Lemmas 2 and 3, we get
\begin{equation*}
\hskip -1.25pc |a_{\lambda}| \leq  C_{\alpha,\beta}
2^{-(\frac{n}{2}+\alpha)j}  2^{-(\frac{n}{2}+|\beta|)j'}
(1+|2^{-j'}k'|) ^{m+\delta\alpha}.
\end{equation*}
\end{enumerate}
{\it Final step.} We consider the case where $\delta>0$ and we
prove that (3.5) implies (3.4). Let $\sigma(x,\xi)
=\sum_{\lambda\in \Lambda_{0,\delta}} a_{\lambda}
\Phi_{\lambda}(x,\xi)$, then we have
\begin{equation*}
\hskip -1.25pc \partial^{\beta}_{\xi} \sigma(x,\xi) = \sum\limits
_{\lambda \in \Lambda_{0,\delta}} 2^{j'|\beta|} a_{\lambda}
\Phi^{\epsilon}_{j,k}(x) (\partial^{\beta}_{\xi}
\Phi^{\epsilon'})_{j',k'}(\xi).
\end{equation*}
By Lemma 2, we have
\begin{align*}
\hskip -1.25pc \|\partial^{\beta}_{\xi}
\sigma_{2}(x,\xi)\|_{B^{\alpha,\infty}_{\infty}} &=
\sup\limits_{\epsilon,j,k} |2^{(\frac{n}{2}+\alpha)j}
\sum\limits_{\epsilon',j',k'} 2^{j'|\beta|} a_{\lambda}
(\partial^{\beta}_{\xi} \Phi^{\epsilon'})_{j',k'} (\xi)|\\[.4pc]
&\leq \sup\limits_{ k} \sum\limits_{\epsilon',j',k'} 2^{j'|\beta|}
|a^{0,\epsilon'}_{0,k,j',k'}| |(\partial^{\beta}_{\xi}
\Phi^{\epsilon'})_{j',k'} (\xi)|\\[.4pc]
&\quad\, + \sup\limits_{\epsilon\neq 0,j,k}
2^{(\frac{n}{2}+\alpha)j} \sum\limits_{\epsilon',j',k'}
2^{j'|\beta|} |a^{\epsilon,\epsilon'}_{j,k,j',k'}|
|(\partial^{\beta}_{\xi} \Phi^{\epsilon'})_{j',k'} (\xi)|\\[.4pc]
&\leq C(1+|\xi|)^{m}+C(1+|\xi|)^{m+\delta\alpha}\leq
C(1+|\xi|)^{m+\delta\alpha} .
\end{align*}
Hence we get $\sigma(x,\xi) \in \tilde{S}^{m}_{0,\delta}$.
\end{proof}

%%pari 11-19

\section{Wavelet characterization for
$\pmb{\tilde{S}^{m}_{\rho,\delta} (\rho > 0)}$}

\setcounter{equation}{0}

We use the wavelet basis which comes from the tensor product of wavelet
basis in $n$-dimension. Let $\Lambda_{\rho,\delta} = \Lambda_{n}\times \Lambda_{n}$. For
$\lambda=(\epsilon,j,k,\epsilon',j',k')\in \Lambda_{\rho,\delta}$,
let $\Phi_{\lambda}(x,\xi) = \Phi^{\epsilon}_{j,k}(x)
\Phi^{\epsilon'}_{j',k'}(\xi)$. Then
$\{\Phi_{\lambda}(x,\xi)\}_{\lambda\in \Lambda_{\rho,\delta}}$ is
an orthogonal normal wavelet basis in $L^{2}(R^{n}\times R^{n})$.
For $a_{\lambda} = \langle
\sigma(x,\xi),\Phi_{\lambda}(x,\xi)\rangle$, it is clear that $\{
a_{\lambda}\}_{\lambda\in \Lambda_{\rho,\delta}}$ becomes a new
representation for symbol. For $\lambda=
(\epsilon,j,k,\epsilon',j',k')$, if $\epsilon'=0$, then $j'=0$ and
we write $\lambda= (\epsilon,j,k,k')$ and $a_{\lambda}=
a^{\epsilon}_{j,k,k'}$. Let $\tau^{\beta}=\tau^{\beta}_{+}$ be the
operator acting on $k'$. We say that $\{ a_{\lambda}\}_{\lambda\in
\Lambda_{\rho,\delta}}\in N^{m}_{\rho,\delta}$, if $a_{\lambda}$
satisfies the following properties:

\begin{enumerate}
\renewcommand\labelenumi{(\roman{enumi})}
\leftskip .1pc
\item The absolute value of $a_{\lambda}$ satisfies:
\begin{equation}
\hskip -4pc |a_{\lambda}| \leq \begin{cases}
C_{\alpha,\beta} 2^{-(\frac{n}{2}+\alpha)j}
2^{-(\frac{n}{2}+\beta)j'} (1+|2^{-j'}k'|)
^{m+\delta\alpha-\rho\beta}, \forall \alpha,\beta\geq 0, &\hbox{if} \
\epsilon' \neq 0; \\[.3pc] C_{\alpha} 2^{-(\frac{n}{2}+\alpha)j}
(1+|k'|) ^{m+\delta\alpha},\forall \alpha\geq 0, & \hbox{if} \
\epsilon'= 0 .
\end{cases}
\end{equation}

\item In addition, for $\epsilon'=0$, $a_{\lambda}$ also satisfies
\begin{equation}
\hskip -4pc |\tau^{\beta}a^{\epsilon}_{j,k,k'}| \leq C_{\alpha,\beta}
2^{-(\frac{n}{2}+\alpha)j} (1+|k'|)
^{m+\delta\alpha-\rho|\beta|},\forall \alpha\geq 0, \beta\in
N^{n}, \ \hbox{if} \ \epsilon'= 0.
\end{equation}
\end{enumerate}

Then we have the following theorem.

\begin{theor}[\!]
The following two conditions are equivalent{\rm :}
\begin{align}
\sigma(x,\xi)\in \tilde{S}^{m}_{\rho,\delta},\\[.3pc]
\{a_{\lambda} \}_{\lambda\in \Lambda_{\rho,\delta}}\in
N^{m}_{\rho,\delta}.
\end{align}
\end{theor}

\begin{proof} {\it From symbol to number array.} To prove (4.1), we
consider first the case where $ \epsilon' \neq 0$. For arbitrary
$\alpha$ and $\beta\in N^{n}$, for sufficiently large
$N'>n+|m|+\delta|\alpha|-\rho|\beta|$, we have
\begin{align}
|a_{\lambda}| &= |\langle \sigma (x,\xi),
\Phi^{\epsilon}_{j,k}(x) \Phi^{\epsilon'}_{j',k'}(\xi)\rangle|\nonumber\\[.3pc]
&= 2^{\frac{nj}{2}} 2^{j'(\frac{n}{2}-|\beta|)} |\langle
\partial^{\beta}_{\xi_{\tau_{\epsilon'}}} \sigma (x,\xi),
\Phi^{\epsilon}(2^{j}x-k) (I^{\beta}_{\epsilon'}
\Phi^{\epsilon'})(2^{j'}\xi-k')\rangle|.
\end{align}
By Lemma~2, we get
\begin{equation*}
|a_{\lambda}| \leq C
2^{-j(\frac{n}{2}+|\alpha|)} 2^{j'(\frac{n}{2}-|\beta|)} \int \frac
{(1+|\xi|)^{m+\delta\alpha-\rho|\beta|}} {(1+|2^{j'}\xi-k'|)^{N'}} \d
\xi.
\end{equation*}
Then applying Lemma~3 to $(1+|\xi|)^{m +\delta\alpha-\rho|\beta|}$, we have
\begin{align*}
|a_{\lambda}| &\leq C
2^{j(\frac{n}{2}-|\alpha|)} 2^{j'(\frac{n}{2}-|\beta|)}
(1+|2^{-j'}k'|)^{m+\delta\alpha-\rho|\beta|}\\[.3pc]
&\quad\, \times \int (1+
|2^{j'}\xi-k'|) ^{|m + \delta\alpha-\rho|\beta||-N'} \d\xi\\[.3pc]
&\leq C 2^{-j(\frac{n}{2}+|\alpha|)} 2^{-j'(\frac{n}{2}+|\beta|)}
(1+|2^{-j'}k'|)^{m+\delta\alpha-\rho|\beta|}.
\end{align*}

For $ \epsilon' = 0$, for arbitrary $\alpha$ and for sufficiently
large $N'>n+|m|+\delta|\alpha|$, by Lemma~2, we have
\begin{align*}
|a_{\lambda}| &= |\langle \sigma (x,\xi),
\Phi^{\epsilon}_{j,k}(x) \Phi^{0}_{0,k'}(\xi)\rangle | \\[.3pc]
&\leq C 2^{-j(\frac{n}{2}+|\alpha|)} \int \frac{(1+|\xi|)^{m+\delta\alpha
-\rho|\beta|}} {(1+| \xi-k'|)^{N'}} \ \d \xi.
\end{align*}
Then applying Lemma~3 to $(1+|\xi|)^{m +\delta\alpha}$, we have
\begin{align*}
|a_{\lambda}| &\leq C
2^{j(\frac{n}{2}-|\alpha|)} (1+| k'|)^{|m+\delta\alpha |}
\int (1+ | \xi-k'|) ^{|m+\delta\alpha |-N'} \ \d\xi\\[.3pc]
&\leq C 2^{-j(\frac{n}{2}+|\alpha|)} (1+| k'|)^{m+\delta\alpha}.
\end{align*}

To prove (4.2), $\forall \alpha\in N^{n}$, let
$\tau^{\alpha}_{\pm} f(x)= \prod^{n}_{i=1}
\tau^{\alpha_{i}}_{\pm e_{i}} f(x)$. Hence we have
\begin{align*}
\tau^{\beta} a^{\epsilon}_{j,k,k'} &= \langle \sigma (x,\xi),
\Phi^{\epsilon}_{j,k}(x) \tau^{\beta} \Phi^{0} (\xi-k')\rangle\\[.3pc]
&= \langle \tau^{\beta}_{-} \sigma (x,\xi), \Phi^{\epsilon}_{j,k}(x)
\Phi^{0} (\xi-k')\rangle.
\end{align*}
For $\xi\in R^{n}, \beta\in N^{n}$,
there exists a $\xi'\in B(\xi, 1+|\beta|)$ such that
$\tau^{\beta}_{-} \sigma (x,\xi) = \partial^{\beta} _{\xi} \sigma
(x,\xi')$. Hence we have
\begin{equation*}
\tau^{\beta} a^{\epsilon}_{j,k,k'} =
\langle \partial^{\beta}_{\xi} \sigma (x,\xi'),
\Phi^{\epsilon}_{j,k}(x) \Phi^{0} (\xi-k')\rangle.
\end{equation*}
Then applying the same argument as above, we get the desired conclusion
(4.2).\vspace{.3pc}

\noindent {\it From wavelet representation to symbol representation.} We
consider three cases: (1)~$|\epsilon||\epsilon'|\neq 0$; (2)
$\epsilon=0,\epsilon'\neq 0$; (3) $\epsilon'= 0$. We calculate the
derivation of the following three symbols:
\begin{align*}
\sigma_{1}(x,\xi) &= \sum\limits_{\substack{\lambda\in
\Lambda_{\rho,\delta}\\[.3pc]
|\epsilon||\epsilon'|\,\neq\,0}}
a_{\lambda} \Phi_{\lambda}(x,\xi),\\[.3pc]
\sigma_{2}(x,\xi) &= \sum\limits_{\substack{\lambda\in
\Lambda_{\rho,\delta}\\ \epsilon=\,0,\epsilon'\,\neq\,0}}
a_{\lambda} \Phi_{\lambda}(x,\xi),\\[.3pc]
\sigma_{3}(x,\xi) &= \sum\limits_{\substack{\lambda\in
\Lambda_{\rho,\delta}\\[.2pc] \epsilon'=\,0}} a_{\lambda}
\Phi_{\lambda}(x,\xi).
\end{align*}

We prove that $\sigma_{1}(x,\xi), \sigma_{2}(x,\xi)\in
S^{m}_{\rho,\delta}$ and $\sigma_{3}(x,\xi)\in
\tilde{S}^{m}_{\rho,\delta}$. As for $\sigma_{1}(x,\xi)$, for
arbitrary $\alpha,\beta\in N^{n}$, we choose $s>|\alpha|,
t>|\beta|, \delta (s-\alpha)\leq \rho (t-\beta)$ and $N'>
n+|m+s\delta-t\rho|$. Then we have
\begin{align*}
|\partial^{\alpha}_{x} \partial ^{\beta}_{\xi} \sigma_{1}(x,\xi)| &\leq \sum\limits_{j,j'\geq 0}
2^{j(|\alpha|-s)} 2^{j'(|\beta|-t)} \sum\limits_{\epsilon,k}
|(\partial ^{\alpha}_{x} \Phi^{\epsilon})(2^{j}x-k)|\\[.3pc]
&\quad\, \times \sum\limits_{\epsilon',k'} (1+|2^{-j'}k'|)^{m+s\delta-t\rho}
|(\partial ^{\beta}_{\xi} \Phi^{\epsilon'})(2^{j'}\xi-k')|\\[.3pc]
&\leq \sum\limits_{j,j'\geq 0} 2^{j(|\alpha|-s)} 2^{j'(|\beta|-t)}
\sum\limits_{k'} \frac{(1+|2^{-j'}k'|)^{m+s\delta-t\rho}}
{(1+|2^{j'}\xi-k'|)^{-N'}}.
\end{align*}
Applying Lemma~3 to $(1+|2^{-j'}k'|)^{m+s\delta-t\rho}$, we have
\begin{align*}
|\partial^{\alpha}_{x} \partial ^{\beta}_{\xi}
\sigma_{1}(x,\xi)| &\leq C_{\alpha,\beta} (1+|\xi|)^{m+\delta |
\alpha|-\rho|\beta|} \sum\limits_{j,j'\geq 0} 2^{j(|\alpha|-s)}
2^{j'(|\beta|-t)}\\[.3pc]
&\leq C_{\alpha,\beta} (1+|\xi|)^{m+\delta |\alpha|-\rho|\beta|}.
\end{align*}

As for $\sigma_{2} (x,\xi)$, for arbitrary $\alpha,\beta\in N^{n}$,
we choose $t>|\beta|$ and $N'> n+|m-t\rho|$. Then we have
\begin{align*}
|\partial^{\alpha}_{x} \partial ^{\beta}_{\xi} \sigma_{2}(x,\xi)| &\leq
\sum\limits_{j'\geq 0}  2^{j'(|\beta|-t)} \sum\limits_{k}
|(\partial ^{\alpha}_{x} \Phi^{0})(x-k)|\\[.3pc]
&\quad\, \times \sum\limits_{\epsilon',k'} (1+|2^{-j'}k'|)^{m-t\rho}
|(\partial^{\beta}_{\xi} \Phi^{\epsilon'})(2^{j'}\xi-k')|\\[.3pc]
&\leq \sum\limits_{j'\geq 0} 2^{j'(|\beta|-t)} \sum\limits_{k'}
\frac{(1+|2^{-j'}k'|)^{m+s\delta-t\rho}}
{(1+|2^{j'}\xi-k'|)^{-N'}}.
\end{align*}
Applying Lemma~3 to $(1+|2^{-j'}k'|)^{m-t\rho}$, we have
\begin{align*}
|\partial^{\alpha}_{x} \partial ^{\beta}_{\xi}
\sigma_{2}(x,\xi)| &\leq C_{\alpha,\beta}
(1+|\xi|)^{m-\rho|\beta|} \sum\limits_{j'\geq 0} 2^{j(|\alpha|-s)}
2^{j'(|\beta|-t)}\\[.3pc]
&\leq C_{\alpha,\beta} (1+|\xi|)^{m-\rho|\beta|}.
\end{align*}
As for $\sigma_{3}(x,\xi)$, for arbitrary $\beta\in N^{n}$, we
have
\begin{align*}
\partial^{\beta}_{\xi} \sigma_{3}(x,\xi) &= \sum\limits_{\epsilon,j,k,k'}
a^{\epsilon}_{j,k,k'} \Phi^{\epsilon}_{j,k}(x) (\partial
^{\beta}_{\xi} \Phi^{0})(\xi-k')\\[.3pc]
&= \sum\limits_{\epsilon,j,k,k'}
\tau^{\beta} a^{\epsilon}_{j,k,k'} \Phi^{\epsilon}_{j,k}(x)
S^{\beta} (\partial ^{\beta}_{\xi} \Phi^{0})(\xi-k').
\end{align*}
For all $\alpha\in N$, we choose $N > n+|m + \delta\alpha-\rho
|\beta||$. Applying Lemmas~1 and 2, we have
\begin{align*}
\| \partial ^{\beta}_{\xi}
\sigma_{3}(x,\xi) \|_{B^{\alpha,\infty}_{\infty}} &\leq
C_{\alpha,\beta} \left\| 2^{j(\frac{n}{2}+\alpha)}\sum\limits_{k'}
\tau^{\beta} a^{\epsilon}_{j,k,k'} S^{\beta} (\partial
^{\beta}_{\xi} \Phi^{0})(\xi-k') \right\|_{\infty}\\[.3pc]
&\leq C \left\|\sum\limits_{k'} (1+|k'|)^{m+\delta\alpha-\rho |\beta|}
(1+|\xi-k'|)^{-N} \right\|\\[.3pc]
&\leq  C_{\alpha,\beta} (1+|\xi|)^{m+\delta |
\alpha|-\rho|\beta|}.
\end{align*}
That is to say, $\sigma_{3}(x,\xi)\in \tilde{S}^{m}_{\rho,\delta}$.\vspace{-.5pc}
\end{proof}

\section{Kernel-distribution}

In this section, we consider the kernel-distribution property of
symbols and prove Theorem~2. By Theorem~7, the kernel-distribution
of the symbol operator $\sigma(x,D)$ can be written as
\begin{align*}
k(x,z) = (2\pi)^{-n} \sum\limits_{(\epsilon,j,k; \epsilon',j',k')\in
\Lambda_{n}} a^{\epsilon,\epsilon'}_{j,k,j',k'}
\Phi^{\epsilon}(2^{j}x-k) \hat{\Phi}^{\epsilon'}(2^{-j'}z)\
\e^{i2^{-j'}k'z},
\end{align*}
where
\begin{equation*}
2^{\frac{n}{2}(j'-j)} a^{\epsilon,\epsilon'}_{j,k,j',k'} \in
N^{m}_{\rho,\delta}.
\end{equation*}
We decompose $k(x,z)$ into three parts:
\begin{align*}
k_{1}(x,z) &= (2\pi)^{-n} \sum\limits_{\substack{\lambda\in
\Lambda_{\rho,\delta}\\[.2pc]
|\epsilon||\epsilon'|\neq 0}}
a^{\epsilon,\epsilon'}_{j,k,j',k'} \Phi^{\epsilon}(2^{j}x-k)
\hat{\Phi}^{\epsilon'}(2^{-j'}z) \ \e^{i2^{-j'}k'z},\\[.3pc]
k_{2}(x,z) &= (2\pi)^{-n} \sum\limits_{\substack{\lambda\in
\Lambda_{\rho,\delta}\\[.2pc] \epsilon=0,\epsilon'\neq 0}}
a^{0,\epsilon'}_{0,k,j',k'} \Phi^{0}(x-k)
\hat{\Phi}^{\epsilon'}(2^{-j'}z) \ \e^{i2^{-j'}k'z},\\[.3pc]
k_{3}(x,z) &= (2\pi)^{-n} \sum\limits_{\substack{\lambda\in
\Lambda_{\rho,\delta}\\[.2pc] \epsilon'= 0}}
a^{\epsilon}_{j,k,k'} \ \e^{ik'z} \Phi^{\epsilon}(2^{j}x-k)
\hat{\Phi}^{0}(z).
\end{align*}
Hence, by Meyer's wavelet property, we know that: (i) if $|z|\leq
\frac{\pi}{3}$, then $k_{1}(x,z)=k_{2}(x,z)=0$; and (ii) if $|z|\geq
\frac{4\pi}{3}$, then $k_{3}(x,z)=0$. Now we prove that
\begin{equation*}
|\partial^{\alpha}_{x}\partial^{\beta}_{z}k_{1}(x,z)| +
|\partial^{\alpha}_{x}\partial^{\beta}_{z}k_{2}(x,z)|\leq
C_{\alpha, \beta, N} (1+|z|)^{-N},\quad \forall N>0
\end{equation*}
and
\begin{equation*}
\|\partial^{\beta}_{z}k_{3}(x,z)\|_{B^{\alpha,\infty}_{\infty}}
\leq C_{\alpha, \beta, N} |z|^{-N},\quad \forall N\geq 0
\end{equation*}
and
\begin{equation*}
n+m+\delta\alpha +\max (1,\rho) \beta < N \rho.
\end{equation*}

First, we consider $k_{1}(x,z)$. For $\alpha,\beta\in N^{n},
\forall s_{1}$ and $t_{1}$, we have
\begin{align*}
I_{1} &= |\partial^{\alpha}_{x}\partial^{\beta}_{z}k(x,z)| \leq
C\sum\limits_{j,j'\geq 0} 2^{(|\alpha|-s_{1})j}
2^{-(n+t_{1})j'}\\[.3pc]
&\quad \times \sum\limits_{k'} (1+|2^{-j'}k'|)^{m+s_{1}\delta-t_{1}\rho+|\beta|}\\[.3pc]
&\quad \times \sum\limits_{\epsilon\neq 0 ,k}
|\partial^{\alpha}_{x}\Phi^{\epsilon}(2^{j}x-k)|
\sum\limits_{\epsilon'\neq 0, |\gamma|\leq |\beta|}
|\partial^{\gamma}_{z} \hat{\Phi}^{\epsilon'}(2^{-j'}z)|.
\end{align*}
By choosing $t_{1}\rho> m+s_{1}\delta +|\beta|+n$, we have
\begin{equation*}
\sum\limits_{k'} (1+|2^{-j'}k'|)^{m+s_{1}\delta-t_{1}\rho+|\beta|}\leq C 2^{nj'}.
\end{equation*}
Note that $\sum_{\epsilon,k} |\partial^{\alpha}_{x}\Phi^{\epsilon}
(2^{j}x-k)|\leq C$; hence we have
\begin{equation*}
I_{1} \leq C\sum\limits_{j,j'\geq 0} 2^{(|\alpha|-s_{1})j} 2^{-t_{1}j'}
\sum\limits_{\epsilon'\neq 0} |\hat{\Phi}^{\epsilon'}(2^{-j'}z)|.
\end{equation*}

Since $\Phi^{\epsilon}(x)(\epsilon\neq 0)$ are Meyer's wavelets,
there exists $0 < M' < M$ such that $\forall \alpha\in N^{n}$,
we have ${\rm supp }\partial^{\alpha}_{z}
\hat{\Phi}^{\epsilon}(z)\subset B(0, 2^{M})\backslash
B(0,2^{M'})$. Hence, there exists at most a finite number $j'$
such that $\sum_{\epsilon'\neq 0}
|\hat{\Phi}^{\epsilon'}(2^{-j'}z)|\neq 0$ and $2^{-j'}\sim
C(1+|z|)^{-1}$. By choosing $s_{1}>|\alpha|$ and $t_{1}> \max
\{\frac{m+s_{1}\delta+|\beta|+n}{\rho}, N\}$, we get
\begin{equation*}
I_{1} \leq C_{\alpha, \beta, N} (1+|z|)^{-N}.
\end{equation*}

Secondly, we consider $k_{2}(x,z)$. For $\alpha,\beta\in N^{n},
s_{2}$ and $t_{2}$, we have
\begin{align*}
I_{2} &= |\partial^{\alpha}_{x}\partial^{\beta}_{z}k(x,z)|
\leq C\sum\limits_{j'\geq 0} 2^{-(n+t_{2})j'} \sum\limits_{k'}
(1+|2^{-j'}k'|)^{m-t_{2}\rho+|\beta|}\\[.3pc]
&\quad \times \sum\limits_{k} |\partial^{\alpha}_{x}\Phi^{0}(x-k)| \sum\limits_{\epsilon'\neq 0,
|\gamma|\leq |\beta|}
|\partial^{\gamma}_{z}\hat{\Phi}^{\epsilon'}(2^{-j'}z)|.
\end{align*}

We choose $t_{2}> \max \{\frac{m+|\beta|+n}{\rho}, N\}$ and
applying the same proof as above, we get
\begin{equation*}
I_{2} \leq C_{\alpha, \beta, N} (1+|z|)^{-N}.
\end{equation*}

Finally, we consider $k_{3}(x,z)$. We know that
\begin{equation*}
\sum\limits_{k'} a^{\epsilon}_{j,k,k'} \ \e^{ik'z} = \prod\limits^{n}_{i=1}
(1 - \e^{iz_{i}})^{-\beta_{i}} \sum\limits_{k'}
(\tau^{\beta}a^{\epsilon}_{j,k,k'}) \ \e^{ik'z}
\end{equation*}
and
\begin{equation*}
\partial^{\gamma}_{z} \sum\limits_{k'} (\tau^{\beta} a^{\epsilon}_{j,k,k'}) \
\e^{ik'z} = C_{\gamma} \sum\limits_{k'} k^{\prime \gamma}(\tau^{\beta} a^{\epsilon}_{j,k,k'}) \
\e^{ik'z}.
\end{equation*}
Hence,
\begin{align*}
\partial^{\gamma}_{z} k_{3}(x,z) &= C \sum\limits_{(\epsilon,j,k)\in
\Lambda_{n}} \Phi^{\epsilon}(2^{j}x-k) \sum\limits_{\gamma_{1}+\gamma_{2}
+\gamma_{3}=\gamma} C^{\gamma_{1},\gamma_{2}}_{\gamma}\\[.3pc]
&\quad \times \left(\partial^{\gamma_{1}}_{z} \prod\limits^{n}_{i=1}(1-
\e^{iz_{i}})^{-\beta_{i}}\right) \left(\sum\limits_{k'} k^{\prime \gamma_{2}}
(\tau^{\beta}a^{\epsilon}_{j,k,k'}) \ \e^{ik'z}\right)
\partial^{\gamma_{3}}_{z}\hat{\Phi}^{0}(z).
\end{align*}
By Lemma~2 and by the estimation of $\tau^{\beta} a^{\epsilon}_{j,k,k'}$, we
choose a convenient $\beta\in N^{n}$ such that $m+\delta\alpha-\rho|\beta|+|
\gamma_{2}|+n<0$ and get the desired conclusion.

\section{$\pmb{L^{2}}$-continuity for symbol operator}

First, we prove a useful lemma. For $i = 1, 2$, let
$\Phi^{i}(x)$ be real-valued functions which belong to $S(R^{n})$.
For $j\geq 0$ and $m=(k,l)\in Z_{2n}$, let
\setcounter{equation}{0}
\begin{equation}
T_{j,m}f(x) = \int \e^{ixy} \Phi^{1}_{j,k}(x) \Phi^{2}_{j,l} (y) f(y) \
\d y = \int K_{j,m}(x,y)f(y) \ \d y,
\end{equation}
where
\begin{equation}
K_{j,m}(x,y) = \e^{ixy}\Phi^{1}_{j,k}(x) \Phi^{2}_{j,l}(y).
\end{equation}
The kernel-distribution of the conjugate operator $T_{j,m}^{*}$ is
\begin{equation}
K_{j,m}^{*}(x,y) = \e^{-ixy}\Phi^{2}_{j,l}(x)
\Phi^{1}_{j,k}(y).
\end{equation}
Then we have the following lemma.

\begin{lem}

$\forall m=(k,l),m'=(k',l')\in Z_{2n}${\rm ,} there exists a
sufficiently large $N_{0}>n$ such that $T_{j,m}$ satisfies the
following two conditions{\rm :}
\begin{align}
\|T_{j,k,l}T^{*}_{j,k',l'} \|_{L^{2}\rightarrow L^{2}} \leq
C(1+4^{-j}|k-k'|)^{-2N_{0}} (1+|l-l'|)^{-2N_{0}},\\[.3pc]
\|T^{*}_{j,k,l}T_{j,k',l'} \|_{L^{2}\rightarrow L^{2}} \leq
C(1+|k-k'|)^{-2N_{0}} (1+4^{-j}|l-l'|)^{-2N_{0}}.
\end{align}
\end{lem}

\vspace{.3pc}
\begin{proof}

The kernel-distribution of $T_{j,k,l}^{*}T_{j,k',l'}$ is
\begin{equation*}
K^{1,j}_{m,m'}(y,z) = \Phi_{1,k,k'}(2^{-j}(y-z))
\Phi^{2}_{j,l}(y)\Phi^{2}_{j,l'}(z),
\end{equation*}
where
\begin{equation*}
\Phi_{1,k,k'}(z) = \int \Phi^{1}(x-k) \Phi^{1}(x-k') \ \e^{-ixz} \ \d x.
\end{equation*}
Since
\begin{equation*}
|\Phi_{1,k,k'}(2^{-j}(y-z))|\leq C(1+|k-k'|)^{N}
(1+2^{-j}|y-z|)^{-N},
\end{equation*}
we get the desired conclusion for the norm of $T_{j,k,l}^{*}T_{j,k',l'}$.

Further, the kernel-distribution of $T_{j,k,l}T_{j,k',l'}^{*}$ is
\begin{equation*}
K^{2,j}_{m,m'}(y,z) = \Phi_{2,l,l'}(2^{-j}(y-z))
\Phi^{1}_{j,k}(y)\Phi^{1}_{j,k'}(z),
\end{equation*}
where
\begin{equation*}
\Phi_{2,l,l'}(z) = \int \Phi^{2}(x-l) \Phi^{2}(x-l') \ \e^{ixz}\  \d x.
\end{equation*}
Since
\begin{equation*}
|\Phi_{2,l,l'}(2^{-j}(y-z))|\leq C(1+|l-l'|)^{N} (1+2^{-j}|y-z|)^{-N},
\end{equation*}
we get the desired conclusion for the norm of $T_{j,k,l}^{*}T_{j,k',l'}$.
\end{proof}

\setcounter{pot}{2}
\begin{pot}
{\rm Let
\begin{equation*}
\tilde{K}^{\epsilon,\epsilon'}_{j}(x,y) =
\e^{ixy}\sum\limits_{k,l} a^{\epsilon,\epsilon'}_{j,k,l}
\Phi^{\epsilon}_{j,k}(x) \Phi^{\epsilon'}_{j,l}(y)
\end{equation*}
be the kernel-distribution of $\tilde{T}^{\epsilon,\epsilon'}_{j}$. By
Lemmas~4 and 5, we have
\begin{equation*}
\|\tilde{T}^{\epsilon,\epsilon'}_{j}\|
_{L^{2}\rightarrow L^{2}}\leq C4^{jn}\sup\limits_{k,l}|a^{\epsilon,\epsilon'}_{j,k,l}|.
\end{equation*}

Let
\begin{equation*}
\tilde{K}(x,y) = \e^{ixy}\sum\limits_{\epsilon,\epsilon',j,k,l}
a^{\epsilon,\epsilon'}_{j,k,l} \Phi^{\epsilon}_{j,k}(x)
\Phi^{\epsilon'}_{j,l}(y)
\end{equation*}
be the kernel-distribution of $\tilde{T}$; then we have
\begin{equation*}
\|\tilde{T}\| _{L^{2}\rightarrow L^{2}}\leq C\sum\limits_{j}4^{jn}
\sup\limits_{\epsilon,\epsilon',k,l} |a^{\epsilon,\epsilon'}_{j,k,l}|.
\end{equation*}

$\left.\right.$\vspace{-1.5pc}

\pagebreak

Let $Ff(x)$ be the Fourier transform of $f(x)$; then we have
$\sigma(x,D)f(x)=\tilde{T} Ff(x)$, i.e. $\sigma(x,D)$ is
continuous from $L^{2}$ to $L^{2}$.

Now we prove part (ii) of Theorem~3. Let $\Phi^{1}(x)$ be a
regular mother wavelet, and ${\rm supp }\Phi^{1}(x)\subset
B(0,2^{M})$ where $M$ is an integer. Let $\tilde{\Phi}(x)$ be the
Fourier transform of the function $(\Phi^{1}(x))^{2}$; then there
exists $C_{1}$ such that for $|x|\leq 2 C_{1}$, we have
$\tilde{\Phi}(x)\geq C_{1}$. For $j\geq 0$, let $\tau_{j}$ be the
set of $l$ satisfying $2^{-M-2}l\in Z^{n}$ and $|l|\leq
C_{2}4^{j}$ where $C_{2}$ satisfies $|\sum_{l\in\tau_{j}}
{\rm e}^{i4^{-j}lx}|\geq C4^{jn}$ for $|x|\leq C_{1}$. Further, for
$j\geq 0, 2^{-M-2} k$ and $2^{-M-2}l\in Z^{n}$, let
$a_{j,k,l} = \e^{-i4^{-j}kl}$; otherwise, $a_{j,k,l} = 0$. For $j \geq
0, l\in \tau_{j}$, let $a_{j,l}=2^{-jn}$; otherwise, $a_{j,l} = 0$.

To show that the result in Theorem~3 is sharp, we construct a
special function and a special operator. Let $f_{j}(x)=\sum_{l}
a_{j,l}\Phi^{1}_{j,l}(x)$ and let
\begin{equation*}
K_{j}(x,y) = \e^{ixy} \sum\limits_{k,l}a_{j,k,l}
\Phi^{1}_{j,k}(x) \Phi^{1}_{j,l}(y)
\end{equation*}
be the kernel-distribution of the operator $\tilde{T}_{j}$. We have
$\|f_{j}\|_{L^{2}}\sim C$ and
\begin{align*}
I_{j} &= \|\tilde{T}_{j}f_{j}(x)\|_{L^{2}}^{2}\\[.3pc]
&= \sum\limits_{k}\int \bigg|\int \sum\limits_{l} a_{j,k,l}a_{j,l}
(\Phi^{1}_{j,l}(y))^{2} \ \e^{ixy} \ \d y\bigg|^{2} (\Phi^{1}_{j,k}(x))^{2}\
\d x.
\end{align*}

According to the definition of $\tilde{\Phi}(x)$ and $a_{j,k,l}$,
we have
\begin{align*}
I_{j} &= \sum\limits_{k} \int \bigg|
\sum\limits_{l} a_{j,k,l}a_{j,l} \ \e^{i2^{-j}lx}\bigg|
\tilde{\Phi}(2^{-j}x)|^{2} (\Phi^{1}_{j,k}(x))^{2}\ \d x\\[.3pc]
&= \sum\limits_{k} \int \bigg| \sum\limits_{l} a_{j,l}\
\e^{i2^{-j}l(x-2^{-j}k)}\bigg|\tilde{\Phi}(2^{-j}x)|^{2}
(\Phi^{1}_{j,k}(x))^{2}\ \d x.
\end{align*}

By changing variables $2^{j}x-k\rightarrow x$ and by the definition of
$a_{j,l}$, we have
\begin{equation*}
I_{j} = 4^{-jn} \sum\limits_{k}\int \bigg|
\sum\limits_{l\in\tau_{j}}\ \e^{i4^{-j}lx}\bigg|^{2} |
\tilde{\Phi}(4^{-j}x+4^{-j}k)|^{2}
|\Phi^{1}(x)|^{2}\ \d x.
\end{equation*}
For $|k|\leq C_{1}4^{j}$ and $|x|\leq
C_{1}$, we have $|\tilde{\Phi}(4^{-j}x+4^{-j}k)|^{2}\geq C_{1}$.
Hence, we have
\begin{equation*}
I_{j} \geq C \int \bigg|\sum\limits_{l\in\tau_{j}}
{\rm e}^{i4^{-j}lx}\bigg|^{2} |\Phi^{1}(x)|^{2} \ \d x\geq C 4^{2jn}.
\end{equation*}

Let $K_{j}(x,\xi)$ be the symbol of the operator $\sigma_{j}(x,D)$;
then we have
\begin{equation*}
\|\sigma_{j}(x,D)\|_{L^{2}\rightarrow L^{2}}\geq C 4^{jn} \quad
\hbox{and} \quad \|K_{j}(x,\xi)\|_{B^{s,\infty}_{\infty}}= 2^{j(s+n)}.
\end{equation*}
That is to say, for $0 < s < n$, there exists a symbol
$\sigma(x,\xi)\in B^{s,\infty}_{\infty}$ but $\sigma(x,D)$ is not
continuous from $L^{2}$ to $L^{2}$.}
\end{pot}

\section{$\pmb{L^{p}}$-continuity}

We begin with a lemma about the characterization of symbol.
\begin{lem}
If $\sigma(x,\xi)$ satisfies condition {\rm (1.6),} then
\setcounter{equation}{0}
\begin{equation}
\sum\limits_{j,\epsilon,\epsilon'} 2^{nj} \sup\limits_{k}
\sum\limits_{l} |a^{\epsilon,\epsilon'}_{j,k,l}|<\infty.
\end{equation}

In addition{\rm ,} for $0<s<n${\rm ,} the following two conditions are
equivalent{\rm :}
\begin{align}
\sum\limits_{j} 2^{j(n+s)}\omega(j) &< \infty,\\[.3pc]
\sum\limits_{j,\epsilon,\epsilon'} 2^{sj} \sup\limits_{k} \sum\limits_{l}
|a^{\epsilon,\epsilon'}_{j,k,l}| &< \infty.
\end{align}
\end{lem}
\begin{proof}
{\it From wavelet representation to symbol}. That is to say,
we prove that (7.3) implies (7.2). For $j\geq 1, e\in I_{2n}$, we
have
\begin{align*}
\sigma_{j,e}(x,\xi) &= \sum\limits_{j'\geq j} \sum\limits_{(\epsilon,
\epsilon',k,l)} a^{\epsilon,\epsilon'}_{j',k,l} \tau^{n}_{2^{-j}e}
\Phi^{\epsilon,\epsilon'}_{j',k,l} (x,\xi)\\[.3pc]
&\quad + \sum\limits_{j'< j} \sum\limits_{(\epsilon,\epsilon',k,l)}
a^{\epsilon,\epsilon'}_{j',k,l} \tau^{n}_{2^{-j} e}
\Phi^{\epsilon,\epsilon'}_{j',k,l} (x,\xi).
\end{align*}
Hence, we have
\begin{align*}
I_{s,e} &= \sum\limits_{j \geq 1} 2^{j(n+s)}\sup\limits_{m\in Z^{n}}
\int_{2^{-j}m+2^{-j}Q} \ \d x \int_{R^{n}} |\sigma_{j,e}(x,\xi)| \ \d \xi\\[.3pc]
&\leq C \sum\limits_{j\geq 1} 2^{j(n+s)}
\sup\limits_{m\in Z^{n}} \int_{2^{-j}m+2^{-j}Q}
\sum\limits_{j'\geq j} \sum\limits_{(\epsilon,\epsilon',k,l)}
|a^{\epsilon,\epsilon'}_{j',k,l}| |\Phi^{\epsilon}( 2^{j'}x-k)| \ \d x\\[.3pc]
&\quad\, + C\sum\limits_{j\geq 1} 2^{j(n+s)}
\sup\limits_{m\in Z^{n}} \int_{2^{-j}m+2^{-j}Q} \sum\limits_{j'< j}
2^{(j'-j)n}\\[.3pc]
&\quad\, \times \sum\limits_{(\epsilon,\epsilon',k,l)}
|a^{\epsilon,\epsilon'}_{j',k,l}| |\Phi^{\epsilon}(2^{j'}x-k)| \ \d x\\[.3pc]
&\leq C\sum\limits_{j\geq 1} \sum\limits_{j'\geq j} 2^{js}
\sup\limits_{k\in Z^{n}} \sup\limits_{\epsilon,\epsilon'}
\sum\limits_{l} |a^{\epsilon,\epsilon'}_{j',k,l}|\\[.3pc]
&\quad\, + C\sum\limits_{j\geq 1} \sum\limits_{j'< j}
2^{(j'-j)n} 2^{js} \sup\limits_{k\in Z^{n}} \sup\limits_{
\epsilon,\epsilon'} \sum\limits_{l}
|a^{\epsilon,\epsilon'}_{j',k,l}|.
\end{align*}
If $0 < s < n$, then
\begin{equation*}
I_{s,e} \leq  C\sum\limits_{j'} 2^{j's} \sup\limits_{\epsilon,\epsilon',k} \sum\limits_{l}
|a^{\epsilon,\epsilon'}_{j',k,l}|.
\end{equation*}
And further, we have
\begin{align*}
I' &= \sup\limits_{m\in Z^{n}}\int_{m + Q} \ \d x
\int_{R^{n}} |\sigma (x,\xi)| \ \d \xi\\[.3pc]
&\leq C \sup\limits_{m\in Z^{n}} \int_{ m+ Q} \sum\limits_{j \geq 0}
\sum\limits_{(\epsilon,\epsilon',k,l)} |a^{\epsilon,\epsilon'}_{j, k, l}|
|\Phi^{\epsilon}( 2^{j }x-k)| \ \d x\\[.3pc]
&\leq C \sum\limits_{j \geq 0} \sup\limits_{k\in Z^{n},\epsilon,\epsilon'}
\sum\limits_{l} |a^{\epsilon,\epsilon'}_{j ,k,l}|\\[.3pc]
&\leq C\sum\limits_{j } 2^{j s} \sup\limits_{\epsilon,\epsilon',k}
\sum\limits_{l} |a^{\epsilon,\epsilon'}_{j ,k,l}|.
\end{align*}
{\it From symbol to wavelet representation}. For
$(\epsilon,\epsilon',j,k,l)\in \Lambda_{2n}$, we have
\begin{equation*}
|a^{\epsilon,\epsilon'}_{j,k,l}| = |\langle \sigma (x,\xi),
\Phi^{\epsilon,\epsilon'}_{j,k,l}(x,\xi) \rangle|.
\end{equation*}
If $|\epsilon|+|\epsilon'|= 0$, then $j=0$ and we have
\begin{align*}
|a^{0,0}_{0,k,l}| &= |\langle \sigma (x,\xi), \Phi ^{0,0}
(x-k, \xi-l) \rangle|\\[.3pc]
&\leq C\sum\limits_{|k-k'|\leq 2^{M}} \int_{k'+ Q} \int_{R^{n}}
|\sigma (x,\xi)| \ \d x \ \d \xi.
\end{align*}
If $|\epsilon|+|\epsilon'|\neq 0$, according to Lemma~1, we have
\begin{align*}
|a^{\epsilon,\epsilon'}_{j,k,l}| &= 2^{jn} |\langle \sigma (x,\xi),
\tau^{n} _{-2^{-1} e_{(\epsilon,\epsilon')}} \tilde{ \Phi}^{\epsilon,
\epsilon'} (2^{j}x-k,2^{j}\xi-l) \rangle|\\[.3pc]
&= 2^{jn} |\langle \sigma_{1+j,(\epsilon,\epsilon')} (x,\xi),
\tilde{\Phi}^{\epsilon,\epsilon'}(2^{j}x-k,2^{j}\xi-l) \rangle|.
\end{align*}
Hence we get
\begin{equation*}
|a^{\epsilon,\epsilon'}_{j,k,l}|\leq C\sum\limits_{|k-k'|\leq
2^{M}} 2^{jn} \int_{2^{-j}k'+2^{-j}Q} \int _{R^{n}}
|\sigma_{1+j,(\epsilon,\epsilon')} (x,\xi)| \ \d x \ \d \xi.
\end{equation*}
So we get the desired conclusion.
\end{proof}

\begin{pot}
{\rm Let
\begin{equation*}
K^{\epsilon,\epsilon'}_{j}(x,y) = \sum\limits_{k,l}
a^{\epsilon,\epsilon'}_{j,k,l} \Phi^{\epsilon}(2^{j}x-k)
\hat{\Phi}^{\epsilon'}(2^{-j}(x-y)) \ \e^{i 2^{-j}l(x-y)}
\end{equation*}
be the kernel-distribution of the operator $T^{\epsilon,\epsilon'}_{j}$.
We have
\begin{equation*}
|K^{\epsilon,\epsilon'}_{j}(x, y)| \leq C \sum\limits_{k}
\sum\limits_{l} |a^{\epsilon,\epsilon'}_{j, k, l}| |\Phi^{\epsilon}
(2^{j}x-k)| |\hat{\Phi}^{\epsilon'} (2^{-j}(x-y))|.
\end{equation*}
That is,
\begin{equation*}
\int |K^{\epsilon,\epsilon'}_{j}
(x,y)|\ \d x \leq C 2^{jn} \sup\limits_{k} \sum\limits_{l}
|a^{\epsilon,\epsilon'}_{j,k,l}|
\end{equation*}
and
\begin{equation*}
\int |K^{\epsilon,\epsilon'}_{j} (x,y)|\ \d y \leq C 2^{jn}
\sup\limits_{k} \sum\limits_{l} |a^{\epsilon,\epsilon'}_{j,k,l}|.
\end{equation*}
Hence, for $1\leq
p\leq\infty$, $T^{\epsilon,\epsilon'}_{j}$ is continuous from
$L^{p}$ to $L^{p}$.

Let
\begin{equation*}
\Gamma= \{(\epsilon,\epsilon',j), \forall k,l\in Z^{n},
(\epsilon,\epsilon',j,k,l)\in \Lambda_{2n}\}.
\end{equation*}
Hence $\sigma(x,D)=\sum_{(\epsilon,\epsilon',j)\in\Gamma}
T^{\epsilon,\epsilon'}_{j}$ is continuous from $L^{p}$ to $L^{p}$
for $1\leq p\leq\infty$.

Then we prove part~(ii) of Theorem~4. Let $M$ be a sufficiently
big integer, let $\Phi^{1}(x)$ be a regular Daubechies' wavelet
with ${\rm supp}\,\Phi^{1}(x)\subset B(0,2^{M})$ and let
$\Phi^{2}(x)$ be Meyer's wavelet. Moreover, let
\begin{equation*}
\sigma_{j}(x,\xi) = \sum\limits_{2^{M+2}k\in Z^{n}}
\Phi^{1}(2^{j}x-k) \Phi ^{2} (2^{j}\xi)
\end{equation*}
and let
\begin{equation*}
\sigma(x,\xi) = \sum\limits_{(2+M)j\in N} j^{2} 2^{-jn} \sigma_{j}(x,\xi).
\end{equation*}
Then $\sigma(x,\xi)$ satisfies conditions (1.8) and (1.9).}
\end{pot}

\section*{Acknowledgements}

This work is supported by the NNSF of China (No.~10001027), the
innovation funds of Wuhan University and the subject construction funds
of Mathematic and Statistic School, Wuhan University.

\end{document}